\documentclass[sn-mathphys-num]{sn-jnl}

\usepackage{graphicx}%
\usepackage{multirow}%
\usepackage{amsmath,amssymb,amsfonts}%
\usepackage{amsthm}%
\usepackage{mathrsfs}%
\usepackage[title]{appendix}%
\usepackage{xcolor}%
\usepackage{textcomp}%
\usepackage{manyfoot}%
\usepackage{booktabs}%
\usepackage{algorithm}%
\usepackage{algorithmicx}%
\usepackage{algpseudocode}%
\usepackage{listings}%
\usepackage{bbm}

\theoremstyle{thmstyleone}%
\newtheorem{theorem}{Theorem}
\theoremstyle{thmstyletwo}%
\newtheorem{remark}{Remark}%
\newtheorem{prop}{Proposition}

\theoremstyle{thmstylethree}%

\newcommand{\mb}[1]{\mathbbm{#1}}
\newcommand{\mc}[1]{\mathcal{#1}}
\newcommand{\Es}{\mathbf{E}}

\begin{document}

\title[One-step corrected projected stochastic gradient descent for statistical estimation]{One-step corrected projected stochastic gradient descent for statistical estimation}

\author*[1]{\fnm{Alexandre} \sur{Brouste}}\email{alexandre.brouste@univ-lemans.fr}

\author[1]{\fnm{Youssef} \sur{Esstafa}}\email{youssef.esstafa@univ-lemans.fr}

\affil[1]{\orgdiv{Laboratoire Manceau de Math\'ematiques}, \orgname{Le Mans Universit\'e}, \country{France}}

\abstract{A generic, fast and asymptotically efficient method for parametric estimation is described. It is based on the projected stochastic gradient descent on the loglikelihood function corrected by a single step of the Fisher scoring algorithm. We show theoretically and by simulations that it is an interesting alternative to the usual stochastic gradient descent with averaging or the adaptative  stochastic gradient descent.}


\maketitle

\section{Introduction}

The stochastic gradient descent \citep{RM51} for finding the root of a given functional is a widely used method in statistical learning. In the parametric estimation setting, this method leads to a (strongly) consistent estimator  but which is not asymptotically efficient in term of converging rate or in term of asymptotic variance depending on the conditions retained. For sublinear functionals, consistency has been shown with probabilistic arguments in  \citep{W52, B54, GH76}  and asymptotic normality in \citep{C54, HL56, S58, F68}. For more general functionals, the ordinary differential equation method has been developed \citep{L77, KC78} with a boundedness assumption of the random sequence. In order to avoid this quite restrictive assumptions, truncated (or projected) stochastic gradient descent has been proposed \citep{C96, KY03}.

The stochastic gradient descent has been improved in two direction to obtain a statistical procedure with optimal asymptotic rate and variance. On the one hand,  the stochastic gradient with averaging has been studied \citep{R88, P90, PJ92}. On the other hand, the adaptative stochastic gradient \citep{LR79, V67} has been suggested.

In this paper, we propose a fast and asymptotically efficient alternative to averaging or adaptivity. It is based on the one-step procedure. 

The one-step procedure was initially considered in  \citep{LeCam} for the estimation of parameters in independent and identically distributed (i.i.d.) samples. In this procedure, an initial guess estimator is proposed which is fast to be computed but not asymptotically efficient.
Then, a single step of the gradient descent method is done on the log-likelihood function in order
to correct the initial estimation and reach asymptotic efficiency. With
some recent developments, the one-step procedure has been successfully generalized to more sophisticated statistical experiments as diffusion processes \citep{GY20, KU15}, ergodic Markov chains~\citep{KM16}, inhomogeneous Poisson and Hawkes
counting processes \citep{BF23,DGK18}, fractional Gaussian and stable noises observed
at high frequency \citep{BM18, BSV20}.

In the following, Section~\ref{sec:notations} is dedicated to notations and known results of convergence rates for stochastic gradient descent (SGD), stochastic gradient descent with averaging (AVSGD), adaptative gradient descent (ADSGD) and maximum likelihood estimation (MLE). The main result on (strong) consistency and asymptotic normality of the one-step procedure in the multidimensional parameter setting is given in Section~\ref{sec:main}. Monte Carlo simulations are done in Section~\ref{sec:sim} to assess the performance of the proposed statistical procedure (OSSGD) in comparison with SGD, AVSGD, ADSGD and MLE in terms of computation time and asymptotic variance for samples of finite size.

\section{Notations}\label{sec:notations}

In our parametric estimation problem, the observation sample is denoted $X^{(n)}=(X_1,\ldots,X_n)$ and is composed of independent and identically distributed random variables.  The probability density $f(\cdot, u)$ (with respect to some $\sigma$-finite measure) of $X_1$ is parametrized by $u \in \Theta \subset \mb{R}^p$ where $\Theta$ is an open set. The true parameter $\vartheta \in \Theta$ is to be estimated.

The estimation problem of the unknown parameter $\vartheta$ can be seen as finding the minimum of an unknown function $G(u) = \Es_\vartheta \left( - \log f (X_1, u) \right)$
or the root of its gradient
\begin{equation}
M(u) = \Es_\vartheta \left( - \nabla_u \log f (X_1, u) \right).
\end{equation}

The standard statistical procedure to estimate the parameter $\vartheta$ is the maximum likelihood estimator (MLE) defined by
\begin{equation}
\widehat{\vartheta}_n = \arg \max_{u \in \Theta} \frac{1}{n} \sum_{i=1}^n \log f (X_i ,u).
\end{equation}
Under regularity assumptions, the MLE is consistent, asymptotically normal
$$ \sqrt{n} \left( \widehat{\vartheta}_n  - \vartheta \right) \Longrightarrow \mc{N}(0, {\cal I}(\vartheta)^{-1}) $$
where ${\cal I}(\vartheta)$ stands for the Fisher information matrix 
\begin{equation}\label{eq:FisherInf}
\mc{I}(\vartheta) = -\Es_\vartheta \left[ \nabla_{u, u}^2 \log f\left(X_1, \vartheta \right)\right]
\end{equation}
and asymptotically efficient in the minimax sense \citep{IK81}. Here $\Longrightarrow$ is the convergence in law as $n\rightarrow \infty$.
But the MLE is generally not in a closed form and its approximation by a classical gradient descent method can be time consuming for large samples. The moment estimator, which is an other generic methodology, when it has closed-form, is generally not asymptotically efficient \citep{IK81}.

Consequently, the Robins-Monro algorithm \citep{RM51} could be considered to find this root. It is defined recursively by
\begin{equation*} \label{eq:SGD} 
\vartheta_{i+1}=\vartheta_i - \gamma_i (-  \nabla_u \log f (X_{i+1}, \vartheta_i)), \quad 0 \leq  i \leq n-1,
\end{equation*} 
where $(\gamma_i)_i$ is the step sequence and $\vartheta_0$ is the initial value (it may be random but square integrable) of  the procedure.

In our estimation problem, the functional $M$ is generally not sublinear and the sequence  $(\vartheta_i)_i$ cannot be considered as bounded in probability. For instance, direct computations when the distribution of $X_1$ is exponential of rate parameter $\vartheta$ give 
$$ M(u)= \frac{1}{\vartheta} - \frac{1}{u}, \quad u >0.$$

Consequently, projected stochastic gradient descent will be considered in the following. Namely,
\begin{equation} \label{eq:projSGD} 
\vartheta_{i+1}=\Pi_K \left[ \vartheta_i - \gamma_i (-  \nabla_u \log f (X_{i+1}, \vartheta_i)) \right], \quad 0 \leq  i \leq n-1,
\end{equation} 
where $\Pi_K$ is the projection onto the constraint set $K = \left\{ u : a_j \leq u_j \leq b_j  \right\}$
for $-\infty < a_j < b_j< \infty$, $j=1,\ldots, p$. It can be reformulated as 
\begin{equation} \label{eq:projSGD2} 
\vartheta_{i+1}= \vartheta_i - \gamma_i (-  \nabla_u \log f (X_{i+1}, \vartheta_i)) +\gamma_i  Z_i , \quad 0 \leq  i \leq n-1,
\end{equation} 
where $\gamma_i  Z_i$ is the shortest Euclidian length to take back $\vartheta_i - \gamma_i (-  \nabla_u \log f (X_{i+1}, \vartheta_i))$ to the constraint set $K$ if it is not in $K$.

Under general assumptions, this procedure leads to a strongly consistent estimator \citep{KY03} for
\begin{equation} \label{eq:stepsize}
\gamma_i \geq 0, \quad \sum_i \gamma_i^2 < \infty \quad \mbox{and} \quad  \sum_i \gamma_i= \infty,
\end{equation}
that is $\vartheta_n \longrightarrow \vartheta$ as $n\rightarrow \infty$ with probability one.
This algorithm is fast but is not asymptotically efficient, neither in terms of converging rate nor in terms of asymptotic variance. For the sequence $\gamma_i= i^{-r}$ and $r \in (1/2,1)$, it leads to an asymptotically normal estimator for which
\begin{equation} \label{eq:TCL1} n^{\frac{r}{2}} \left( \vartheta_n - \vartheta \right)  \Longrightarrow \mc{N}\left(0, \frac{1}{2}I_p\right)
\end{equation}
where $I_p$ stands for the $p\times p$ identity matrix. It is worth mentioning that the asymptotic variance does not depend on $\vartheta$ in the i.i.d. setting.  

It had also been shown that the stochastic gradient descent with $\gamma_i= c i^{-1}$ and
\begin{equation}\label{eq:const}
c> \frac{1}{2\lambda_{min}(\mc{I}(\vartheta))},
\end{equation}
where $\lambda_{min}(A)$ is the lowest eigenvalue of the matrix $A$, is asymptotically rate efficient but is still not asymptotically variance efficient with
$$ \sqrt{n}\left( \vartheta_n - \vartheta \right)  \Longrightarrow \mc{N}\left(0, c^2 \mc{I}(\vartheta) (2c \mc{I}(\vartheta) - I_p)^{-1}\right).
$$
The constraint \eqref{eq:const} depends on the unknown parameter and cannot be used in practice.

Consequently, in order to speed up the estimation convergence rate in \eqref{eq:TCL1}, two methods are classically used: averaging and adaptivity.

\subsection{Averaging}

The averaging method was proposed (see \cite{PJ92,R88,KY03}) to reach variance efficiency with
$$ \overline{\vartheta}_n = \frac{1}{n} \sum_{i=1}^n \vartheta_i.$$
This estimator is consistent, asymptotically normal with efficient rate and variance, namely 
$$  \sqrt{n} \left( \overline{\vartheta}_n - \vartheta \right) \Longrightarrow \mc{N}(0, {\cal I(\vartheta)}^{-1}).$$
 
\subsection{Adaptivity}

 In the simple setting of i.i.d. samples, the adaptative stochastic gradient descent writes
$$\widetilde{\vartheta}_{i+1}=\widetilde{\vartheta}_i - i^{-1}{\cal I}(\widetilde{\vartheta}_i)^{-1} \left( -  \nabla_u  \log f( X_{i+1}, \widetilde{\vartheta}_i) \right), \quad 0 \leq  i \leq n-1.$$
When the classical assumptions are fulfilled, it leads also to a consistent and asymptotical normal estimators with optimal limit variance (see \cite{A98} and the references therein), namely
$$  \sqrt{n} \left( \widetilde{\vartheta}_n - \vartheta \right) \Longrightarrow \mc{N}(0, {\cal I(\vartheta)}^{-1}).$$

\section{One-step correction}\label{sec:main}

In order to improve the convergence rate of the gradient descent algorithm, we propose in the following the one-step procedure starting from an initial guess estimator taken from the projected stochastic gradient algorithm. This procedure is shown to be faster than the classical computation of the MLE but still asymptotically efficient. It is an interesting alternative to the stochastic gradient algorithm with averaging or adaptative gradient descent and shows nice properties also on samples of finite size.

In the one-step estimation procedure, the estimation $\vartheta_n$ given at step $n$ by the projected stochastic gradient descent \eqref{eq:SGD} is corrected by
\begin{equation}\label{eq:OS} 
\vartheta^*_n = \vartheta_n + \mc{I}(\vartheta_n)^{-1} \cdot \frac{1}{n}\sum_{i=1}^n \nabla_u \log f (X_i, \vartheta_n).
\end{equation}
It leads to a consistent, asymptotically normal and asymptotically efficient estimator of $\vartheta$ (see Theorem 1 below). 

In the following, we recall the slow convergence of the projected stochastic gradient descent in the multidimensional setting which is the initial guess estimator in the one-step procedure.

Let $Y_i= \nabla_u  \log f( X_{i+1}, \vartheta_i)$ and $\Es_j$ the conditional expectation with respect to the $\sigma$-algebra generated by $\{\vartheta_0, (Y_i, i < j) \}$.  Let $\gamma_i= i^{-r}$ and $r \in (1/2,1)$ in the algorithm~\eqref{eq:projSGD}. The classical assumptions are formulated in \cite[Section 10.4 p. 341]{KY03}, namely

\begin{enumerate}

\item[A.1] The true value $\vartheta$ is in the interior of the constraint set $K$ and $\vartheta_n \rightarrow \vartheta$ as $n\rightarrow \infty$ with probability one;

\smallskip

\item[A.2] For small $\rho>0$, $\left\{ Y_n \mathbb{I}_{\{ |\vartheta_n - \vartheta| \leq \rho \}} \right\}$ is uniformly integrable and there is a function $g$ such that for $ |\vartheta_n - \vartheta| \leq \rho$,
$$   \Es_n Y_n = g(\vartheta_n); $$

\smallskip

\item[A.3] There exists a constant  $0 < C < \infty$ such that for small $\rho >0$, 
$$ \sup_n \Es_n |Y_n|^2 \mathbb{I}_{\{ |\vartheta_n - \vartheta| \leq \rho \}} < C \quad \mbox{w.p.1;} $$

\smallskip

\item[A.4]There is a Hurwitz matrix $A$ such that 
$$ g(u) = A (u-\vartheta) + o(|u-\vartheta |).$$ 
\end{enumerate}

The following proposition gives the $n^{\frac{r}{2}} $-consistency of the initial guess estimator which is the first key ingredient in order to prove our next Theorem~\ref{thm:OS}.

\smallskip

\begin{prop}[\cite{KY03}]\label{prop:SGD}  Under aforementioned assumptions, the sequence $ n^{\frac{r}{2}} \left( \vartheta_n - \vartheta \right)$ is tight.
\end{prop}

\smallskip

\begin{remark}
Additive assumptions in order to fullfill A.1 are given in \cite[Section 5.2 p 125]{KY03}.
\end{remark}

\smallskip

\begin{remark}
Additive  assumptions to obtain the asymptotic normality,
\begin{equation} n^{\frac{r}{2}} \left( \vartheta_n - \vartheta \right)  \Longrightarrow \mc{N}\left(0, \frac12 I_p\right),
\end{equation}
are given in \cite[Section 10.2 p. 329]{KY03}. 
\end{remark}

This algorithm is fast but is not asymptotically efficient, neither in terms of converging rate nor in terms of asymptotic variance. In order to obtain asymptotic normality with optimal rate and variance for the one-step corrected projected stochastic gradient descent,  we also suppose that
\begin{enumerate}
\item[A.5.] The matrix valued function ${\cal I}(\vartheta)$ is Lipschitz continuous, {\it i.e.} there exists a constant $L>0$ such that
$$ \left\| {\cal I}(x)-{\cal I}(y) \right\|_m \leq  L \| x - y \|   , \quad x, y \in \Theta,$$
where $\| \cdot \|_m$ and $\| \cdot \|$ stand for Euclidean norms in the space of matrices and vectors respectively. 
\end{enumerate}

With this condition, we can state the main result:

\smallskip

\begin{theorem}\label{thm:OS} The sequence  $(\vartheta^*_n, n \geq 1)$ of one-step estimators of $\vartheta$ defined by~\eqref{eq:OS} is consistent and asympotically normal, {\it i.e.}
\begin{equation}
\sqrt{n} ( \vartheta^*_n -  \vartheta) \Longrightarrow \mc{N}(0, \mc{I}(\vartheta)^{-1}).
\end{equation}
\end{theorem}

It is worth emphasizing that both speed and asymptotic variance are improved in this one-step procedure due to the regularity of the Fisher information matrix.

\begin{proof}
The proof is postponed in Appendix~\ref{app:OS}.
\end{proof}

\section{Simulations}\label{sec:sim}

The joint estimation of the shape parameter $\alpha$  and scale parameter $\beta$ is considered in the statistical experiment generated by a sample $X^{(n)}=(X_1,X_2,\ldots,X_n)$ of i.i.d. Gamma random variables whose probability density function is given by
$$ f(x)= \frac{\beta^\alpha}{\Gamma(\alpha)}x^{\alpha-1} \exp( - \beta x ), \quad x >0.$$

Let us denote $\vartheta=(\alpha,\beta)$. In this statistical experiment, the sequence of maximum likelihood estimators $(\widehat{\vartheta}_n)_{n \geq 1}$ of $\vartheta$ is not in a closed-form.  The sequence of MLE satisfies 
$$ \sqrt{n}\left( \widehat{\vartheta}_n-\vartheta \right) \rightarrow {\cal N}\left(0,{\cal I}(\vartheta)^{-1}\right),$$
where 
$$ {\cal I}(\vartheta)=\begin{pmatrix} \psi^{(2)}(\alpha) & -\frac{1}{\beta} \\ - \frac{1}{\beta} &  \frac{\alpha}{\beta^2}\end{pmatrix}.$$
Here, $\psi^{(n)}$ is the polygamma functions (see  \cite[Section 6.4.1, page 260]{AS72}) defined by $\psi^{(n)}(\alpha)=\frac{\partial^n}{\partial \alpha^n} \log \Gamma(\alpha)$.

The different estimators (MLE, SGD, OSSGD, AVSGD, ADSGD) have been compared in terms of variance and computation time on $B=2\times 10^3$ Monte Carlo simulations for samples of size $n=2\times 10^4$. The SGD is done with $\gamma_i=i^{-r}$ where $r$ is chosen to be equal to $0.6$. It is worth mentioning that the results are similar for all values of $\frac 12 < r < 1$.

\begin{figure}
\centering
\includegraphics[width=\textwidth]{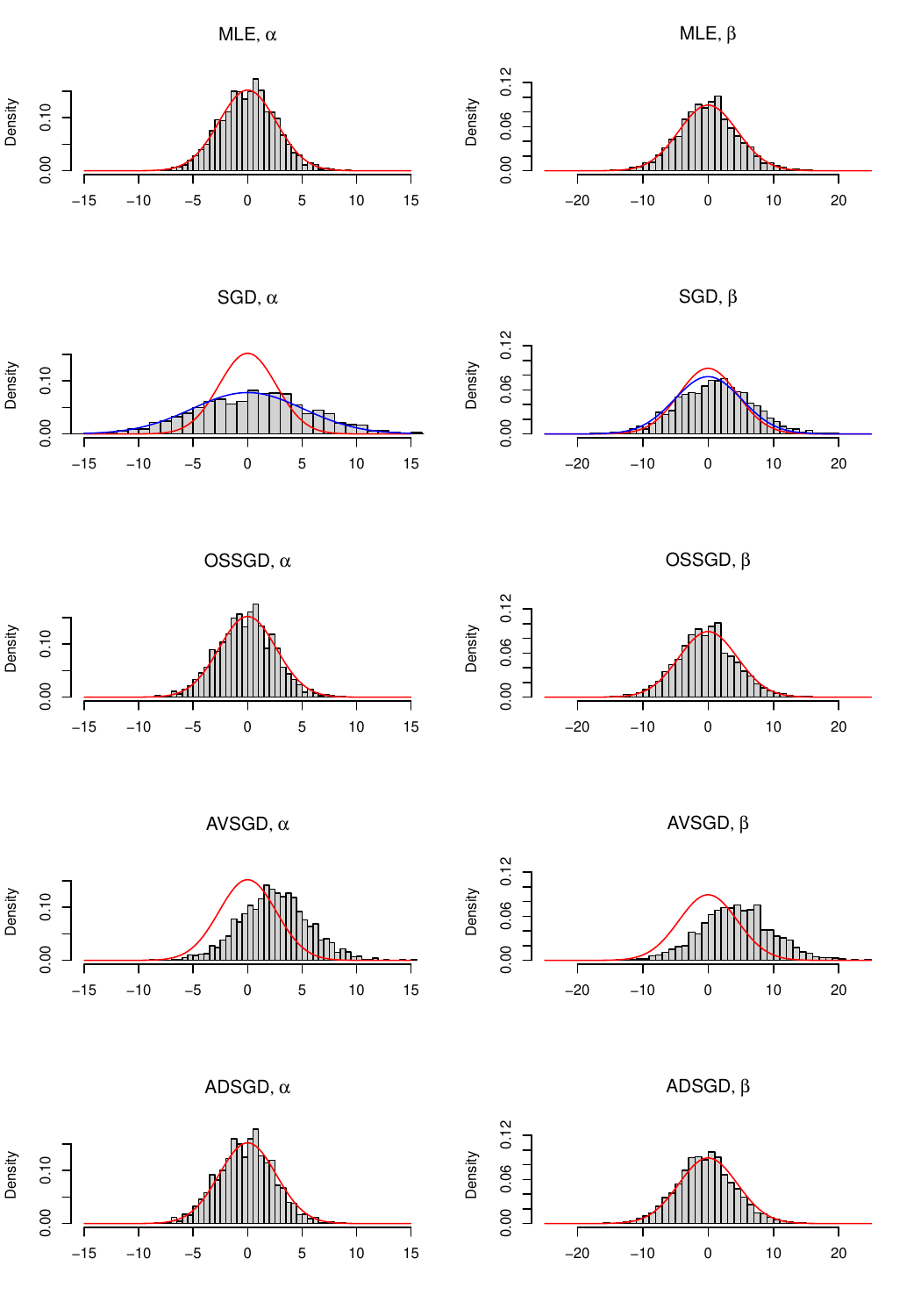}
\caption{Statistical errors renormalized by $\sqrt{n}$ for MLE, SGD, OSSGD, AVSGD and ADSGD for $n=2\times10^4$ and $B=2\times 10^3$ Monte-Carlo simulations. Theoretical optimal variance (in red) and variance of SGD (in blue) are superimposed. }\label{fig:hist2}
\end{figure}

We can see on Figure~\ref{fig:hist2} that the optimal variance is reached for the OSSGD (as for the MLE, AVSGD and ADSGD) that naturally overperforms the non-optimal variance of the slowly converging SGD.  It is worth noting the relative bias for samples of finite size of the AVSGD when the initial value $\vartheta_0$ is fixed.

In terms of computation time, the OSSGD (as the AVSGD) is more than 3 times faster than the MLE. In comparison, the ADSGD is more than two times faster. 

For these reasons, the fast and asymptotically efficient OSSGD is a proper alternative to the averaged and the adapted stochastic gradient descent methods.


\begin{table}[ht]
\centering
\begin{tabular}{rrrrrr}
  \hline
 & MLE & SGD & OSSGD & AVSGD & ADSGD \\ 
  \hline
time (s) & 198.22 & 63.87 & 64.25 & 64.22 & 88.20 \\ 
   \hline
\end{tabular}
\end{table}

It can also be noticed that, for the specific case of the estimation of the parameters in the Gamma distribution, moment estimators \citep{BDN21} or other original explicit estimators \citep{YC17} could have been considered as initial guess estimation in the one-step procedure instead of the SGD. 

\section{Conclusion}

In this paper, we propose to apply the one-step procedure to the slowly converging stochastic gradient descent in order to improve the convergence rate and reach asymptotical efficiency. It is a fast and asymptotically efficient alternative to averaging or adaptivity.

The one-step procedure for the stochastic gradient descent is considered here in the i.i.d. setting but it will be extended in a further work to the regression setting (linear regression, logistic regression (see also \cite{BGP20} for an adaptative procedure), generalized linear models) for larger applications.




\appendix

\section{One-step procedure}\label{app:OS}

For an observation sample $(X_1, \ldots, X_n)$, let us denote $\ell_n(u)= \sum_{i=1}^n \log f (X_i, u)$. Recall that $\vartheta$ is the true parameter and 
\begin{equation}\label{LC1bis}
\vartheta^*_n  =   \vartheta_n  +  {\cal I}(\vartheta_n)^{-1} \cdot  \frac{1}{n} \nabla_u \ell_n(\vartheta_n), \quad n \geq 1.
 \end{equation}

\paragraph{Consistency:}  The consistency of the sequence of initial guess estimators gives, as $n\rightarrow \infty$,  $\vartheta_n \longrightarrow \vartheta $ in probability.  Since $\Es_\vartheta \nabla_u \ell_n(\vartheta) =  0$, the uniform law of large numbers gives, as $n\rightarrow \infty$,
$$  \frac1n \nabla_u \ell_n(\vartheta_n) \longrightarrow 0_{\mathbb{R}^p}
$$
in probability. The uniform continuity of the Fisher information matrix gives the result. Since the initial stochastic gradient descent is also strongly consistent~\citep{KY03}, we can also obtain the strong consistency with the strong law of large numbers.

\paragraph{Asymptotic normality:}
From \eqref{LC1bis}, we have 
$$ \sqrt{n} \left( \vartheta^*_n - \vartheta \right) =  \sqrt{n} \left( \vartheta_n - \vartheta \right) + {\cal I}(\vartheta_n)^{-1} \cdot  \frac{1}{\sqrt{n}} \nabla_u \ell_n(\vartheta_n). $$
The mean-value theorem gives
$$  \nabla_u \ell_n(\vartheta_n) =  \nabla_u \ell_n(\vartheta) + \int_0^1   \nabla_{u,u}^2 \ell_n  \left( \vartheta + \tau ( \vartheta_n - \vartheta) \right) d\tau \cdot \left( \vartheta_n - \vartheta \right)
 $$ 
and
\begin{align} \label{eq:LC} \sqrt{n} \left( \vartheta^*_n - \vartheta \right) &=   n^\frac{r}{2}  \left\lbrace  I_p + {\cal I}(\vartheta_n)^{-1} \frac{1}{n} \int_0^1   \nabla_{u,u}^2 \ell_n  \left( \vartheta + \tau ( \vartheta_n - \vartheta) \right) d\tau \right\rbrace  n^\frac{r}{2} \left( \vartheta_n - \vartheta \right)n ^{\frac12-r} \nonumber   \\
&\quad +  {\cal I}(\vartheta_n)^{-1} \cdot  \frac{1}{\sqrt{n}} \nabla_u \ell_n(\vartheta),
\end{align}
where $I_p$ is the $p\times p$ identity matrix.\\
The central limit theorem gives, as $n\rightarrow \infty$,
$$ \frac{1}{\sqrt{n}} \nabla_u \ell_n(\vartheta)  \Longrightarrow {\cal N} \left( 0,  {\cal I}(\vartheta) \right)$$
in law and the proper convergence of the second term in the r.h.s. of Equation~\eqref{eq:LC}.

Considering the first right-hand term, we have that $(\vartheta_n)_{n \geq 1}$ is $n^\frac{r}{2}$-consistent by assumption and $n ^{\frac12-r} \rightarrow 0$, as $n\rightarrow \infty$, for $\frac12 < r \leq 1$. Then, we need to show that
$$n^\frac{r}{2}  \left(  I_p + {\cal I}(\vartheta_n)^{-1}  \frac{1}{n} \int_0^1   \nabla_{u,u}^2 \ell_n  \left( \vartheta + \tau ( \vartheta_n - \vartheta) \right) d\tau\right) = n^\frac{r}{2}  A $$ 
is bounded in probability as $n\rightarrow \infty$ with 
\begin{eqnarray*}
A & =& {\cal I}(\vartheta_n)^{-1} \left( {\cal I}(\vartheta_n) +\frac{1}{n}\int_0^1   \nabla_{u,u}^2 \ell_n  \left( \vartheta + \tau ( \vartheta_n - \vartheta) \right) d\tau \right) \nonumber \\
&=& {\cal I}(\vartheta_n)^{-1} \cdot \left(  \left[  {\cal I}(\vartheta_n)- {\cal I}(\vartheta)  \right] +  \left[   \frac1n \nabla_{u,u}^2  \ell_n(\vartheta)   +{\cal I}(\vartheta) \right]  \right. \nonumber\\
&&  \left. +  \frac1n\int_0^1   \left[  \frac{\partial^2}{\partial \vartheta^2} \ell_n  \left( \vartheta + \tau (\vartheta_n - \vartheta) \right)  -  \nabla_{u,u}^2 \ell_n(\vartheta)  \right] d\tau  \right).\end{eqnarray*}

The second terms in the r.h.s. converges  to zero at rate $\sqrt{n}$. The Lipschitz continuity of the Fisher information allows to control the first and third terms by $C \| \vartheta_n  - \vartheta \|$ where $C$ is a generic constant. Since $(\vartheta_n)_{n \geq 1}$ is $n^\frac{r}{2}$-consistent, the quantity $n^\frac{r}{2} A$ is bounded in probability as $n\rightarrow \infty$. The Slutsky theorem gives the final result.

\paragraph{Acknowledgments}
We would like to thank Alain Bensoussan for all the fruitful discussions on this work.
This research partially benefited from the support of the ANR pro\-ject 'Efficient inference for large and high-frequency data' (ANR-21-CE40-0021) and the ’Efficience et Sobriété Numériques’ research program under the aegis of Fondation du Risque, a joint initiative by Le Mans University and EREN Groupe.

\bibliography{Biblio.bib}

\end{document}